\documentclass{amsart}
\usepackage{amsfonts,amsmath,amsthm,amssymb}
\newtheorem{theorem}{Theorem}
\newtheorem{corollary}{Corollary}
\newtheorem{proposition}{Proposition}
\title[Torsion units]
{Torsion units in integral group ring of the \\ Mathieu simple group $M_{22}$}
\date{}

\author{V.A.~Bovdi, A.B.~Konovalov and S.~Linton}

\address{V.A.~Bovdi
\newline Institute of Mathematics, University of Debrecen,
\newline P.O.  Box 12, H-4010 Debrecen, Hungary
\newline Institute of Mathematics and Informatics, College of Ny\'\i regyh\'aza,
\newline S\'ost\'oi \'ut 31/b, H-4410 Ny\'\i regyh\'aza, Hungary}
\email{vbovdi@math.klte.hu}

\address{A.B.~Konovalov
\newline School of Computer Science,
         University of St Andrews,
\newline North Haugh, St Andrews, Fife, KY16 9SX, Scotland}
\email{konovalov@member.ams.org}

\address{S.~Linton
\newline School of Computer Science,
         University of St Andrews,
\newline North Haugh, St Andrews, Fife, KY16 9SX, Scotland}
\email{sal@cs.st-and.ac.uk}

\thanks{The research was supported by OTKA grants No.T 43034,
No.K61007 and Francqui Stichting (Belgium) grant ADSI107}



\begin{document}
\begin{abstract}
We investigate the possible character values of 
torsion units of the normalized unit group of the
integral group ring of Mathieu sporadic group
$M_{22}$. We confirm the Kimmerle conjecture on prime graphs for this
group and specify the partial augmentations for possible
counterexamples to the stronger Zassenhaus conjecture.
\end{abstract}

\maketitle

\section{Introduction, conjectures   and main results}
\label{Intro}

Let $V(\mathbb Z G)$ be  the normalized unit group of the
integral group ring $\mathbb Z G$ of  a finite group $G$. \, A long-standing
conjecture  of H.~Zassenhaus \, {\bf (ZC)} \, says that every torsion unit
$u \in V(\mathbb ZG)$ is conjugate within the rational group algebra
$\mathbb Q G$ to an element in $G$ (see \cite{Zassenhaus}).

For finite simple groups the main tool for the investigation of the
Zassenhaus conjecture is the Luthar--Passi method, introduced in
\cite{Luthar-Passi} for the case of $A_{5}$ and then applied in 
\cite{Luthar-Trama} for the case of $S_{5}$.
Later M.~Hertweck in \cite{Hertweck1} extended
the  method and applied it to $PSL(2,p^{n})$.
The same method has also proved to be useful for some groups containing 
non-trivial normal subgroups. For some recent results we refer to
\cite{Bovdi-Hofert-Kimmerle,Bovdi-Konovalov,Hertweck2,Hertweck3,
Hertweck1,Hofert-Kimmerle}. 
Some related properties and   weakened variations of the
Zassenhaus  conjecture can be found in
\cite{Artamonov-Bovdi, Bleher-Kimmerle,Kimmerle}.

To define the conjectures we will investigate, and describe the
methods we will use,  we  introduce some notation. By $\# (G)$ we
denote the set of all primes dividing the order of $G$. The
Gruenberg--Kegel graph (or the prime graph) of $G$ is the graph
$\pi (G)$ with vertices labeled by  $\# (G)$ and 
an edge from $p$ to $q$ if there is an element of order $pq$ in
 $G$. In \cite{Kimmerle} W.~Kimmerle   proposed the
following weakened variation of the Zassenhaus conjecture:

\begin{itemize}
\item[]{\bf (KC)} \qquad
If $G$ is a finite group then $\pi (G) =\pi (V(\mathbb Z G))$.
\end{itemize}

It is easy to see that {\bf (ZC)} implies {\bf{KC}} since it implies
that the set of orders of torsion units of $V(\mathbb Z G)$ is the
same the set of orders of elements of $G$.

In the same paper W.~Kimmerle verified that {\bf (KC)} holds for
finite Frobenius and solvable groups.
We remark that with respect to the
so-called $p$-version of the Zassenhaus conjecture the investigation
of Frobenius groups was completed by M.~Hertweck and the first author
in \cite{Bovdi-Hertweck}. 
In \cite{Bovdi-Jespers-Konovalov,Bovdi-Konovalov,
Bovdi-Konovalov-M23,Bovdi-Konovalov-Siciliano}
{\bf (KC)} was confirmed for the sporadic simple groups $M_{11}$,
$M_{12}$, $M_{23}$ and some Janko simple groups.

Here we continue these investigations for the Mathieu simple group
$M_{22}$. Although we cannot prove the rational conjugacy of torsion
units of $V(\mathbb Z M_{22})$ with elements of $M_{22}$, our main
result gives a lot of information on the orders and partial 
augmentations of these units. In particular, we confirm  Kimmerle's 
conjecture for this group.

Let $G=M_{22}$. It is well known
(see \cite{GAP,Witt}) that 
$|G| = 2^7 \cdot 3^2 \cdot 5 \cdot 7 \cdot 11$ 
and
$exp(G) = 2^3 \cdot 3 \cdot  5 \cdot 7 \cdot 11$. 
The group $G$ has 12 irreducible characters of the 
following degrees: 
$1$, $21$, $45$, $45$, $55$, $99$, $154$, $210$, 
$231$, $280$, $280$ and $385$.
Let 
$$
\mathcal{C} =\{ C_{1}, C_{2a}, C_{3a}, C_{4a}, C_{4b}, C_{5a},
C_{6a},C_{7a}, C_{7b}, C_{8a}, C_{11a},C_{11b} \}
$$
be the collection of all conjugacy classes of $M_{22}$, where the first
index denotes the order of the elements of this conjugacy class
and $C_{1}=\{ 1\}$. Suppose $u=\sum \alpha_g g \in V(\mathbb Z G)$
has finite order $k$. Denote by
$\nu_{nt}=\nu_{nt}(u)=\varepsilon_{C_{nt}}(u)=\sum_{g\in C_{nt}}
\alpha_{g}$ the partial augmentation of $u$ with respect to
$C_{nt}$. From the Berman--Higman Theorem (see \cite{Berman} and \cite{Sandling}, Ch.5, p.102)  one knows that
$\nu_1 =\alpha_{1}=0$ and
\begin{equation}\label{E:1}
\sum_{C_{nt}\in \mathcal{C}} \nu_{nt}=1.
\end{equation}
Hence, for any character $\chi$ of $G$, we get that $\chi(u)=\sum
\nu_{nt}\chi(h_{nt})$, where $h_{nt}$ is a representative of the
conjugacy class $ C_{nt}$.

Our main result is the following
\begin{theorem}\label{T:1}
Let $G$ denote the Mathieu simple group $M_{22}$. Let  $u$ be
a torsion unit of $V(\mathbb ZG)$ of order $|u|$. Denote by
$\frak{P}(u)$ the tuple
$$
(\nu_{2a},\nu_{3a},\nu_{4a},\nu_{4b},\nu_{5a},\nu_{6a},\nu_{7a},\nu_{7b},\nu_{8a},\nu_{11a},\nu_{11b})\in \mathbb Z^{11}
$$
of partial augmentations of $u$ in $V(\mathbb ZG)$. 
The following properties hold.

\begin{itemize}

\item[(i)] There are no elements of order 
$10$, $14$, $15$, $21$, $22$, $33$, $35$, 
$55$ or $77$ in $V(\mathbb ZG)$. Equivalently, 
if $|u| \not \in \{12,24\}$, then $|u|$ is the 
order of some element $g \in G$. 

\item[(ii)] If $|u| \in \{2,3,5\}$, then $u$ is rationally conjugate to
some $g\in G$.

\item[(iii)] If $|u|=4$, then all components of $\frak{P}(u)$ are zero
except possibly $\nu_{2a}$,$\nu_{4a}$ and $\nu_{4b}$, and the triple
$(\nu_{2a},\nu_{4a},\nu_{4b})$ is one of 
\[
\begin{split}
 \{ \; &
( -2, -1, 4 ), \; ( 2, -1, 0 ), \; ( 0, -1, 2 ), \; ( 0, 6, -5 ), \;
( -2, 6, -3 ), \; ( 0, 5, -4 ), \; \\ 
& ( -2, 5, -2 ),  ( -2, 4, -1 ),  ( 2, 4, -5 ),  ( 0, 4, -3 ), 
  ( -2, 3, 0 ),  ( 2, 3, -4 ),  ( 0, 3, -2 ), \\
& ( 2, 0, -1 ), \; ( -2, 0, 3 ), \; ( 0, 0, 1 ), \; ( 0, -6, 7 ), \; 
  ( 2, -6, 5 ), \; ( 0, 2, -1 ), \; ( 2, 2, -3 ), \\
& ( -2, 2, 1 ), ( 0, -5, 6 ), \; ( 2, -5, 4 ), \; ( 0, -4, 5 ), \;
  ( 2, -4, 3 ), \; ( -2, -3, 6 ), \; ( 0, -3, 4 ), \\
& ( 2, -3, 2 ), \; ( -2, -2, 5 ), \; ( 0, -2, 3 ), \; ( 2, -2, 1 ), \; 
  ( 0, 1, 0 ), \; ( -2, 1, 2 ), \; ( 2, 1, -2 ) \; \} .
\end{split}
\]

\item[(iv)] If $|u|=6$, then all components of $\frak{P}(u)$ are zero
except possibly $\nu_{2a}$,$\nu_{3a}$ and $\nu_{6a}$, and the triple
$(\nu_{2a},\nu_{3a},\nu_{6a})$ is one of 
\[
\begin{split}
 \{ \; 
& ( -4, 6, -1 ), \; ( -2, 6, -3 ), \; ( 4, -9, 6 ), \; ( -4, 9, -4 ), \; ( -2, 3, 0 ), \\
& ( -4, 3, 2 ), \; ( 0, 3, -2 ), \; ( 2, 0, -1 ), \; ( -2, 0, 3 ), \;  ( 0, 0, 1 ), \\
& ( 2, -6, 5 ), \; ( 4, -6, 3 ), \; ( 0, -3, 4 ), \; ( 4, -3, 0 ), \; ( 2, -3, 2 ) \; \}.
\end{split}
\]

\item[(v)] If $|u|=7$, then all components of $\frak{P}(u)$ are zero
except possibly $\nu_{7a}$ and $\nu_{7b}$ and the pair $(\nu_{7a},\nu_{7b})$
is one of 
\[
\begin{split}
  \{ \; (0,1), \; (2,-1), \; (1,0), & \; (-1,2) \; \}.
\end{split}
\]

\item[(vi)] If $|u|=11$, then all components of $\frak{P}(u)$ are zero
except possibly $\nu_{11a}$ and $\nu_{11b}$ and the pair $(\nu_{11a},\nu_{11b})$
is one of 
\[
\begin{split}
 \{ \; & ( 5, -4 ), \, ( 0, 1 ), \, ( -2, 3 ), \, ( 2, -1 ), \, ( -3, 4 ), \\
       & ( -4, 5 ), \, ( 1, 0 ), \, ( 3, -2 ), \, ( -1, 2 ), \, ( 4, -3 ) \; \}.
\end{split}
\]

\end{itemize}

\end{theorem}

{\bf Note} that using our implementation of the Luthar--Passi method,
which we intend to make available in the GAP package LAGUNA \cite{LAGUNA}, 
we are able to compute the set of 76 tuples containing (likely as a proper
subset) possible tuples of partial augmentations for units of order 8, 
listed in the Appendix 1. For the case of order 12 in the Appendix 2 we 
listed 1166 tuples which can not be eliminated using the Luthar--Passi method.

As an immediate consequence of  part (i) of the Theorem we obtain

\begin{corollary} 
If $G \cong M_{22}$ then $\pi(G)=\pi(V(\mathbb ZG))$.
\end{corollary}

\section{Preliminaries}
The following result allows a reformulation of
the Zassenhaus conjecture in terms of vanishing of
partial augmentations of torsion units.

\begin{proposition}\label{P:5}
(see \cite{Luthar-Passi} and
Theorem 2.5 in \cite{Marciniak-Ritter-Sehgal-Weiss})
Let $u\in V(\mathbb Z G)$
be of order $k$. Then $u$ is conjugate in $\mathbb
QG$ to an element $g \in G$ if and only if for
each $d$ dividing $k$ there is precisely one
conjugacy class $C$ with partial augmentation
$\varepsilon_{C}(u^d) \neq 0 $.
\end{proposition}

The next results now serve to restrict the possible values of the  partial augmentations
of torsion units.

\begin{proposition}\label{P:4}
(see \cite{Hertweck2}, Proposition 3.1;
\cite{Hertweck1}, Proposition 2.2)
Let $G$ be a finite
group and let $u$ be a torsion unit in $V(\mathbb
ZG)$. If $x$ is an element of $G$ whose $p$-part,
for some prime $p$, has order strictly greater
than the order of the $p$-part of $u$, then
$\varepsilon_x(u)=0$.
\end{proposition}

The next result is explained in detail in \cite{Luthar-Passi} and
\cite{Bovdi-Hertweck,Hertweck1}.
\begin{proposition}\label{P:1}
(see \cite{Hertweck1,Luthar-Passi}) Let either $p=0$ or $p$ a prime
divisor of $|G|$. Suppose
that $u\in V( \mathbb Z G) $ has finite order $k$ and assume $k$ and
$p$ are coprime in case $p\neq 0$. If $z$ is a complex primitive $k$-th root
of unity and $\chi$ is either a classical character or a $p$-Brauer
character of $G$, then for every integer $l$ the number
\begin{equation}\label{E:2}
\mu_l(u,\chi, p ) =
\textstyle\frac{1}{k} \sum_{d|k}Tr_{ \mathbb Q (z^d)/ \mathbb Q }
\{\chi(u^d)z^{-dl}\}
\end{equation}
is a non-negative integer.
\end{proposition}

Note that if $p=0$, we will use the notation $\mu_l(u,\chi,*)$ 
for $\mu_l(u,\chi , 0)$.

Finally, we shall use the well-known bound for
orders of torsion units.

\begin{proposition}\label{P:2}  (see  \cite{Cohn-Livingstone})
The order of a torsion element $u\in V(\mathbb ZG)$
is a divisor of the exponent of $G$.
\end{proposition}

In case of units of prime power order, the following Proposition
may also be useful to eliminate some tuples of partial augmentations.

\begin{proposition}\label{P:6} (see \cite{Cohn-Livingstone})
Let $p$ be a prime, and let $u$ be a torsion unit of $V(\mathbb ZG)$ of order $p^n$.
Then for $m \ne n$ the sum of all partial augmentations of $u$ with respect to
conjugacy classes of elements of order $p^m$ is divisible by $p$.
\end{proposition}

\section{Proof of the Theorem}

Throughout this section we denote $M_{22}$ by $G$. 
The ordinary and  $p$-Brauer character tables of $G$, which will be denoted by
$\mathfrak{BCT}{(p)}$ where $p\in\{2,3,5,7,11\}$, can be found using
the computational algebra system GAP \cite{GAP}, which derives its
data from \cite{AFG,ABC}. For the characters
and conjugacy classes we will use throughout the paper the same
notation,  including indexation, as used in the GAP Character
Table Library.

Since the group $G$ possesses elements of
orders $2$, $3$, $4$, $5$, $6$, $7$, $8$ and  $11$, we first 
investigate  units of all of these orders except 8.
After this, since by Proposition \ref{P:2}, the
order of each torsion unit divides the exponent of $G$, it
remains to consider units of orders $10$, $12$, $14$, $15$, $21$, 
$22$, $33$, $35$, $55$ and $77$. We prove that $V(\mathbb ZG)$
contains no units of any of these 
orders, except possibly for orders $12$ and $24$.

Now we consider each case separately.

\noindent $\bullet$ Let $u$ be a unit of order $2$, $3$ or $5$. Using
Proposition \ref{P:4} we immediately obtain that all partial
augmentations except one are zero. Thus by Proposition \ref{P:5} part
(ii) of  Theorem \ref{T:1} is proved.

\noindent $\bullet$ Let $u$ be a unit of order $4$. 
By (\ref{E:1}) and Proposition \ref{P:4} we get
$ \nu_{2a}+\nu_{4a}+\nu_{4b}=1$. Now using Proposition \ref{P:1} 
we obtain the following system of inequalities:
\[
\begin{split}
\mu_{0}(u,\chi_{2},*) & = \textstyle \frac{1}{4} (10 \nu_{2a} + 2 \nu_{4a} + 2 \nu_{4b} + 26) \geq 0; \\
\mu_{2}(u,\chi_{2},*) & = \textstyle \frac{1}{4} (-10 \nu_{2a} - 2 \nu_{4a} - 2 \nu_{4b} + 26) \geq 0; \\
\mu_{0}(u,\chi_{5},*) & = \textstyle \frac{1}{4} (14 \nu_{2a} + 6 \nu_{4a} - 2 \nu_{4b} + 62) \geq 0; \\
\mu_{2}(u,\chi_{5},*) & = \textstyle \frac{1}{4} (-14 \nu_{2a} - 6 \nu_{4a} + 2 \nu_{4b} + 62) \geq 0. \\
\end{split}
\]
Put $t_1 = 5 \nu_{2a} +  \nu_{4a} +  \nu_{4b}$ and $t_2 = 7 \nu_{2a} + 3 \nu_{4a} -  \nu_{4b}$,
then $t_1 \in \{ 2r+1 \mid -7 \le r \le 6 \}$ and
$t_2 \in \{ 2s+1 \mid -16 \le s \le 15 \}$. Thus, we obtain the system of linear equations 
$
\nu_{2a}+\nu_{4a}+\nu_{4b} =1, \quad
5 \nu_{2a} +  \nu_{4a} +  \nu_{4b}  = t_1, \quad
7 \nu_{2a} + 3 \nu_{4a} -  \nu_{4b} = t_2.
$
Solving such systems for all possible combinations of values of $t_1$ and $t_2$, 
and considering additional inequalities
\[
\begin{split}
\mu_{0}(u,\chi_{5},3) & = \textstyle \frac{1}{4} (2 \nu_{2a} - 6 \nu_{4a} + 2 \nu_{4b} + 50) \geq 0; \\
\mu_{2}(u,\chi_{5},3) & = \textstyle \frac{1}{4} (-2 \nu_{2a} + 6 \nu_{4a} - 2 \nu_{4b} + 50) \geq 0, \\
\end{split}
\]
and also restrictions given by Proposition \ref{P:6}, we get only
the 34 integer solutions $(\nu_{2a},\nu_{4a},\nu_{4b})$ listed in part (iii) of the Theorem \ref{T:1}.

\noindent $\bullet$ Let $u$ be a unit of order $6$.
By (\ref{E:1}) and Proposition \ref{P:4} we get
$ \nu_{2a}+\nu_{3a}+\nu_{6a}=1$. By Proposition \ref{P:1}
we obtain the following system of inequalities:
\[
\begin{split}
\mu_{1}(u,\chi_{2},*) & = \textstyle \frac{1}{6} (5 \nu_{2a} + 3 \nu_{3a} -  \nu_{6a} + 13) \geq 0; \\
\mu_{3}(u,\chi_{2},*) & = \textstyle \frac{1}{6} (-10 \nu_{2a} - 6 \nu_{3a} + 2 \nu_{6a} + 22) \geq 0; \\
\mu_{0}(u,\chi_{4},7) & = \textstyle \frac{1}{6} (12 \nu_{2a} + 60) \geq 0; \quad 
\mu_{3}(u,\chi_{4},7)   = \textstyle \frac{1}{6} (-12 \nu_{2a} + 48) \geq 0; \\
\mu_{1}(u,\chi_{3},*) & = \textstyle \frac{1}{6} (-3 \nu_{2a} + 48) \geq 0. \\
\end{split}
\]
Using calculations similar to the previous case,
we get only the 15 integer solutions $(\nu_{2a},\nu_{3a},\nu_{6a})$ listed in part (iv) 
of the Theorem \ref{T:1}.

\noindent $\bullet$ Let $u$ be a unit of order $7$.
By (\ref{E:1}) and Proposition \ref{P:4} we get
$\nu_{7a}+\nu_{7b}=1$. Using Proposition \ref{P:1}
we obtain the following system of inequalities:
\[
\begin{split}
\mu_{1}(u,\chi_{3},*) & = \textstyle \frac{1}{7} (4 \nu_{7a} - 3 \nu_{7b} + 45) \geq 0; \quad 
\mu_{3}(u,\chi_{3},*)   = \textstyle \frac{1}{7} (-3 \nu_{7a} + 4 \nu_{7b} + 45) \geq 0; \\
\mu_{1}(u,\chi_{2},2) & = \textstyle \frac{1}{7} (4 \nu_{7a} - 3 \nu_{7b} + 10) \geq 0; \quad 
\mu_{3}(u,\chi_{2},2)   = \textstyle \frac{1}{7} (-3 \nu_{7a} + 4 \nu_{7b} + 10) \geq 0. \\
\end{split}
\]
Using that $\nu_{7a}+\nu_{7b}=1$, we get that $-1 \le \nu_{7a} \le 2$, and after this it is
easy to check that we have only the four integer solutions $(\nu_{7a},\nu_{7b})$ listed in part (v) 
of the Theorem \ref{T:1}.

\noindent $\bullet$ Let $u$ be a unit of order $11$.
By (\ref{E:1}) and Proposition \ref{P:4} we get
$ \nu_{11a}+\nu_{11b}=1$. Using Proposition \ref{P:1}
we obtain the following system of inequalities:
\[
\begin{split}
\mu_{1}(u,\chi_{10},*) & = \textstyle \frac{1}{11} (6 \nu_{11a} - 5 \nu_{11b} + 280) \geq 0; \\
\mu_{2}(u,\chi_{10},*) & = \textstyle \frac{1}{11} (-5 \nu_{11a} + 6 \nu_{11b} + 280) \geq 0; \\
\mu_{1}(u,\chi_{5},2) & = \textstyle \frac{1}{11} (7 \nu_{11a} - 4 \nu_{11b} + 70) \geq 0; \\
\mu_{2}(u,\chi_{5},2) & = \textstyle \frac{1}{11} (-4 \nu_{11a} + 7 \nu_{11b} + 70) \geq 0; \\
\mu_{1}(u,\chi_{5},3) & = \textstyle \frac{1}{11} (6 \nu_{11a} - 5 \nu_{11b} + 49) \geq 0; \\
\mu_{2}(u,\chi_{5},3) & = \textstyle \frac{1}{11} (-5 \nu_{11a} + 6 \nu_{11b} + 49) \geq 0. \\
\end{split}
\]
Using calculations similar to the previous case, we get
only the ten integer solutions for $(\nu_{11a},\nu_{11b})$ listed in part (vi) 
of the Theorem \ref{T:1}.

It remains to prove part (i) of the Theorem \ref{T:1}, 
considering units of $V(\mathbb ZG)$ of orders 
$10$, $14$, $15$, $21$, $22$, $33$, $35$, $55$ and $77$.

\noindent $\bullet$ Let $u$ be a unit of order $10$.
By (\ref{E:1}) and Proposition \ref{P:4} we get
$ \nu_{2a}+\nu_{5a}=1$. Using Proposition \ref{P:1}
we obtain the following system of inequalities:
\[
\begin{split}
\mu_{0}(u,\chi_{2},*) & = \textstyle \frac{1}{10} (20 \nu_{2a} + 4 \nu_{5a} + 30) \geq 0; \\
\mu_{5}(u,\chi_{2},*) & = \textstyle \frac{1}{10} (-20 \nu_{2a} - 4 \nu_{5a} + 20) \geq 0; \\
\mu_{1}(u,\chi_{3},*) & = \textstyle \frac{1}{10} (-3 \nu_{2a} + 48) \geq 0, \\
\end{split}
\]
that has no solutions such that all $\mu_{i}(u,\chi_{j},*)$ are non-negative integers.

\noindent $\bullet$ Let $u$ be a unit of order $14$.
Then by (\ref{E:1}) and Proposition \ref{P:4} we get
$\nu_{2a}+\nu_{7a}+\nu_{7b}=1$.
We need to consider four cases
defined by part (v) of  Theorem \ref{T:1}, but in all of them
using Proposition \ref{P:1} we obtain the same system of inequalities:
\[
\begin{split}
\mu_{0}(u,\chi_{2},*) & = \textstyle \frac{1}{14} (30 \nu_{2a} + 26) \geq 0; \quad 
\mu_{7}(u,\chi_{2},*)   = \textstyle \frac{1}{14} (-30 \nu_{2a} + 16) \geq 0, \\
\end{split}
\]
which has no solutions such that all $\mu_{i}(u,\chi_{j},*)$ are non-negative integers.

\noindent $\bullet$ Let $u$ be a unit of order $15$.
By (\ref{E:1}) and Proposition \ref{P:4} we get
$ \nu_{3a}+\nu_{5a}=1$. Using Proposition \ref{P:1}
we obtain the following system of inequalities:
\[
\begin{split}
\mu_{0}(u,\chi_{2},*) & = \textstyle \frac{1}{15} (24 \nu_{3a} + 8 \nu_{5a} + 31) \geq 0; \\ 
\mu_{5}(u,\chi_{2},*) & = \textstyle \frac{1}{15} (-12 \nu_{3a} - 4 \nu_{5a} + 22) \geq 0. \\
\end{split}
\]
From this follows that $3 \nu_{3a} +  \nu_{5a} = -2$, 
and all conditions together leave us no integer solutions.

\noindent $\bullet$ Let $u$ be a unit of order $21$. Then
by (\ref{E:1}) and Proposition \ref{P:4} we have
$$
\nu_{3a}+\nu_{7a}+\nu_{7b}=1.
$$
We need to consider four cases determined by part (v) of  Theorem \ref{T:1}. 
Using Proposition \ref{P:1} we obtain the following systems of inequalities:
\[
\begin{split}
\mu_{0}(u,\chi_{2},*) & = \textstyle \frac{1}{21} (36 \nu_{3a} + 27) \geq 0; \qquad \quad 
\mu_{7}(u,\chi_{2},*)   = \textstyle \frac{1}{21} (-18 \nu_{3a} + 18) \geq 0; \\
\mu_{1}(u,\chi_{3},*) & = \textstyle \frac{1}{21} (3 \nu_{7a} - 4 \nu_{7b} + \alpha) \geq 0; \quad 
\mu_{9}(u,\chi_{3},*)   = \textstyle \frac{1}{21} (- 6 \nu_{7a} + 8 \nu_{7b} + \alpha) \geq 0; \\
\mu_{3}(u,\chi_{3},*) & = \textstyle \frac{1}{21} (8 \nu_{7a} - 6 \nu_{7b} + \beta) \geq 0, \\
\end{split}
\]
$$
\text{where} \quad 
{\tiny
(\alpha, \beta) = 
\begin{cases}
(49,42), \quad \text{when} \quad \chi(u^3)=\chi(7a); \\
(42,49), \quad \text{when} \quad \chi(u^3)=\chi(7b); \\
(56,35), \quad \text{when} \quad \chi(u^3)= 2\chi(7a)-\chi(7b); \\
(35,56), \quad \text{when} \quad \chi(u^3)= -\chi(7a)+2\chi(7b), \\
\end{cases}}
$$
which have no solutions such that all $\mu_{i}(u,\chi_{j},*)$ are non-negative integers.

\noindent $\bullet$ Let $u$ be a unit of order $22$.
Then by (\ref{E:1}) and Proposition \ref{P:4} we get
$$
\nu_{2a}+\nu_{11a}+\nu_{11b}=1.
$$
We need to consider ten cases determined by part (vi) of the Theorem \ref{T:1}. 
In each case using Proposition \ref{P:1} we obtain the following system of inequalities:
\[
\begin{split}
\mu_{0}(u,\chi_{2},*) & = \textstyle \frac{1}{22} (50 \nu_{2a} - 10 \nu_{11a} - 10 \nu_{11b} + 16) \geq 0; \\
\mu_{11}(u,\chi_{2},*) & = \textstyle \frac{1}{22} (-50 \nu_{2a} + 10 \nu_{11a} + 10 \nu_{11b} + 6) \geq 0, \\
\end{split}
\]
that has no solutions such that all $\mu_{i}(u,\chi_{j},*)$ are non-negative integers.

\noindent $\bullet$ Let $u$ be a unit of order $33$.
Then by (\ref{E:1}) and Proposition \ref{P:4} we get
$$
\nu_{3a}+\nu_{11a}+\nu_{11b}=1.
$$
As in the previous case, we need to consider ten cases determined by part (vi) of the Theorem \ref{T:1}. 
In each case using Proposition \ref{P:1} we obtain the same system:
\[
\begin{split}
\mu_{0}(u,\chi_{2},*) & = \textstyle \frac{1}{33} (60 \nu_{3a} - 20 \nu_{11a} - 20 \nu_{11b} + 17) \geq 0; \\
\mu_{11}(u,\chi_{2},*) & = \textstyle \frac{1}{33} (-30 \nu_{3a} + 10 \nu_{11a} + 10 \nu_{11b} + 8) \geq 0, \\
\end{split}
\]
that has no solutions such that all $\mu_{i}(u,\chi_{j},*)$ are non-negative integers.

\noindent $\bullet$ Let $u$ be a unit of order $35$.
Then by (\ref{E:1}) and Proposition \ref{P:4} we get
$\nu_{5a}+\nu_{7a}+\nu_{7b}=1$.
We need to consider four cases
defined by part (v) of the Theorem \ref{T:1}, but in all of them
using Proposition \ref{P:1} we obtain the same system of inequalities:
\[
\begin{split}
\mu_{0}(u,\chi_{2},*) & = \textstyle \frac{1}{35} (24 \nu_{5a} + 25) \geq 0; \quad 
\mu_{0}(u,\chi_{7},2)   = \textstyle \frac{1}{35} (-48 \nu_{5a} + 90) \geq 0, \\
\end{split}
\]
which has no solutions such that all $\mu_{i}(u,\chi_{j},p)$ are non-negative integers.

\noindent $\bullet$ Let $u$ be a unit of order $55$.
Then by (\ref{E:1}) and Proposition \ref{P:4} we get
$$
\nu_{5a}+\nu_{11a}+\nu_{11b}=1.
$$
We need to consider ten cases determined by part (vi) of the Theorem \ref{T:1}. 
In each case using Proposition \ref{P:1} we obtain the following system of inequalities:
\[
\begin{split}
\mu_{0}(u,\chi_{2},*) & = \textstyle \frac{1}{55} (40 \nu_{5a} - 40 \nu_{11a} - 40 \nu_{11b} + 15) \geq 0; \\
\mu_{11}(u,\chi_{2},*) & = \textstyle \frac{1}{55} (-10 \nu_{5a} + 10 \nu_{11a} + 10 \nu_{11b} + 10) \geq 0; \\
\mu_{1}(u,\chi_{10},*) & = \textstyle \frac{1}{55} (- 6 \nu_{11a} + 5 \nu_{11b} + \alpha) \geq 0; \\
\mu_{5}(u,\chi_{10},*) & = \textstyle \frac{1}{55} (24 \nu_{11a} - 20 \nu_{11b} + \alpha) \geq 0; \\
\mu_{1}(u,\chi_{2},*) & = \textstyle \frac{1}{55} ( \nu_{5a} -  \nu_{11a} -  \nu_{11b} + 21) \geq 0, \\
\end{split}
\]
$$
\text{where} \quad 
{\tiny
\alpha = 
\begin{cases}
286, \quad \text{when} \quad \chi(u^5) = \chi(11a); \\
275, \quad \text{when} \quad \chi(u^5) = \chi(11b); \\
330, \quad \text{when} \quad \chi(u^5) = 5 \chi(11a) - 4 \chi(11b); \\
253, \quad \text{when} \quad \chi(u^5) = -2 \chi(11a) + 3 \chi(11b); \\
297, \quad \text{when} \quad \chi(u^5) = 2 \chi(11a) - \chi(11b); \\
242, \quad \text{when} \quad \chi(u^5) = -3 \chi(11a) + 4 \chi(11b); \\
231, \quad \text{when} \quad \chi(u^5) = -4 \chi(11a) + 5 \chi(11b); \\
308, \quad \text{when} \quad \chi(u^5) = 3 \chi(11a) - 2 \chi(11b); \\
264, \quad \text{when} \quad \chi(u^5) = - \chi(11a) + 2 \chi(11b); \\
319, \quad \text{when} \quad \chi(u^5) = 4 \chi(11a) - 3 \chi(11b), \\
\end{cases}}
$$
that has no solutions such that all $\mu_{i}(u,\chi_{j},*)$ are non-negative integers.

\noindent $\bullet$ Let $u$ be a unit of order $77$. 
Then by (\ref{E:1}) and Proposition \ref{P:4} we have
$$
\nu_{7a}+\nu_{7b}+\nu_{11a}+\nu_{11b}=1.
$$
We must consider $40$ cases determined by parts (v) and (vi) of the Theorem \ref{T:1},
but luckily in all of them using Proposition \ref{P:1} we obtain the same system of inequalities:
\[
\begin{split}
\mu_{11}(u,\chi_{2},*) & = \textstyle \frac{1}{77} (10 \nu_{11a} + 10 \nu_{11b} + 11) \geq 0; \\
\mu_{0}(u,\chi_{2},*) & = \textstyle \frac{1}{77} (- 60 \nu_{11a} - 60 \nu_{11b} + 11) \geq 0, \\
\end{split}
\]
which has no solutions such that all $\mu_{i}(u,\chi_{j},*)$ are non-negative integers.
This finishes the proof of Theorem \ref{T:1}.

\bibliographystyle{plain}
\bibliography{Preprint_Bovdi_Konovalov_Linton_M22}

\vspace{5pt}

\centerline{\bf Appendix 1.}

\vspace{5pt}

Possible partial augmentations $(\nu_{2a},\nu_{4a},\nu_{4b},\nu_{8a})$ for units of order 8:
$$
{\tiny
\begin{array}{lllll}
( -2, -2, 4, 1 ),  & ( -2, -1, 3, 1 ),  & ( -2, 0, 2, 1 ),  & ( -2, 0, 4, -1 ),  & ( -2, 1, 1, 1 ), \\
( -2, 1, 3, -1 ),  & ( -2, 2, 0, 1 ),  & ( -2, 2, 2, -1 ),  & ( -2, 2, 4, -3 ),  & ( -2, 3, -1, 1 ), \\
( -2, 3, 1, -1 ),  & ( -2, 3, 3, -3 ),  & ( -2, 4, -2, 1 ),  & ( -2, 4, 0, -1 ),  & ( -2, 4, 2, -3 ), \\
( -2, 5, -3, 1 ),  & ( -2, 5, -1, -1 ),  & ( -2, 5, 1, -3 ),  & ( -2, 6, -2, -1 ),  & ( -2, 6, 0, -3 ), \\
( 0, -5, 3, 3 ),  & ( 0, -4, 2, 3 ),  & ( 0, -4, 4, 1 ),  & ( 0, -3, 1, 3 ),  & ( 0, -3, 3, 1 ), \\
( 0, -3, 5, -1 ),  & ( 0, -2, 0, 3 ),  & ( 0, -2, 2, 1 ),  & ( 0, -2, 4, -1 ),  & ( 0, -2, 6, -3 ), \\
( 0, -1, -1, 3 ),  & ( 0, -1, 1, 1 ),  & ( 0, -1, 3, -1 ),  & ( 0, -1, 5, -3 ),  & ( 0, 0, -2, 3 ), \\
( 0, 0, 0, 1 ),  & ( 0, 0, 2, -1 ),  & ( 0, 0, 4, -3 ),  & ( 0, 1, -3, 3 ),  & ( 0, 1, -1, 1 ), \\
( 0, 1, 1, -1 ),  & ( 0, 1, 3, -3 ),  & ( 0, 2, -4, 3 ),  & ( 0, 2, -2, 1 ),  & ( 0, 2, 0, -1 ), \\
( 0, 2, 2, -3 ),  & ( 0, 3, -5, 3 ),  & ( 0, 3, -3, 1 ),  & ( 0, 3, -1, -1 ),  & ( 0, 3, 1, -3 ), \\
( 0, 4, -4, 1 ),  & ( 0, 4, -2, -1 ),  & ( 0, 4, 0, -3 ),  & ( 0, 5, -3, -1 ),  & ( 0, 5, -1, -3 ), \\
( 0, 6, -2, -3 ),  & ( 2, -5, 1, 3 ),  & ( 2, -5, 3, 1 ),  & ( 2, -4, 0, 3 ),  & ( 2, -4, 2, 1 ), \\
( 2, -4, 4, -1 ),  & ( 2, -3, -1, 3 ),  & ( 2, -3, 1, 1 ),  & ( 2, -3, 3, -1 ),  & ( 2, -2, -2, 3 ), \\
( 2, -2, 0, 1 ),  & ( 2, -2, 2, -1 ),  & ( 2, -1, -3, 3 ),  & ( 2, -1, -1, 1 ),  & ( 2, -1, 1, -1 ), \\
( 2, 0, -2, 1 ),  & ( 2, 0, 0, -1 ),  & ( 2, 1, -3, 1 ),  & ( 2, 1, -1, -1 ),  & ( 2, 2, -2, -1 ), \\
( 2, 3, -3, -1 ) \\
\end{array}}
$$

\vspace{5pt}

\centerline{\bf Appendix 2.}

\vspace{5pt}

Possible partial augmentations $(\nu_{2a},\nu_{3a},\nu_{4a},\nu_{4b},\nu_{6a})$ for units of order 12:
$$
{\tiny
\begin{array}{llll}
( -4, 9, -4, -4, 4 ),  & ( -4, 9, -3, -5, 4 ),  & ( -4, 9, -2, -6, 4 ),  & ( -4, 9, -1, -7, 4 ), \\
( -4, 9, -1, -3, 0 ),  & ( -4, 9, 0, -8, 4 ),  & ( -4, 9, 0, -4, 0 ),  & ( -4, 9, 1, -9, 4 ), \\
( -4, 9, 2, -10, 4 ),  & ( -4, 9, 3, -11, 4 ),  & ( -4, 12, -5, -5, 3 ),  & ( -4, 12, -4, -6, 3 ), \\
( -4, 12, -4, -4, 1 ),  & ( -4, 12, -3, -7, 3 ),  & ( -4, 12, -3, -5, 1 ),  & ( -4, 12, -2, -8, 3 ), \\
( -4, 12, -2, -6, 1 ),  & ( -4, 12, -1, -9, 3 ),  & ( -4, 12, -1, -7, 1 ),  & ( -4, 12, 0, -10, 3 ), \\
( -4, 12, 0, -8, 1 ),  & ( -3, 3, 0, 1, 0 ),  & ( -3, 3, 1, 0, 0 ),  & ( -3, 3, 2, -1, 0 ), \\
( -3, 3, 2, 3, -4 ),  & ( -3, 3, 3, -2, 0 ),  & ( -3, 3, 3, 2, -4 ),  & ( -3, 3, 4, -3, 0 ), \\
( -3, 3, 4, 1, -4 ),  & ( -3, 3, 5, -4, 0 ),  & ( -3, 3, 5, 0, -4 ),  & ( -3, 3, 6, -5, 0 ), \\
( -3, 6, -4, -1, 3 ),  & ( -3, 6, -3, -2, 3 ),  & ( -3, 6, -2, -3, 3 ),  & ( -3, 6, -2, 1, -1 ), \\
( -3, 6, -1, -4, 3 ),  & ( -3, 6, -1, 0, -1 ),  & ( -3, 6, 0, -5, 3 ),  & ( -3, 6, 0, -1, -1 ), \\
( -3, 6, 1, -6, 3 ),  & ( -3, 6, 1, -2, -1 ),  & ( -3, 6, 2, -7, 3 ),  & ( -3, 6, 2, -3, -1 ), \\
( -3, 6, 3, -8, 3 ),  & ( -3, 6, 3, -4, -1 ),  & ( -3, 6, 4, -9, 3 ),  & ( -3, 6, 4, -5, -1 ), \\
( -3, 6, 5, -10, 3 ),  & ( -3, 6, 5, -6, -1 ),  & ( -3, 6, 6, -11, 3 ),  & ( -3, 9, -7, -4, 6 ), \\
( -3, 9, -6, -5, 6 ),  & ( -3, 9, -6, -3, 4 ),  & ( -3, 9, -5, -6, 6 ),  & ( -3, 9, -5, -4, 4 ), \\
( -3, 9, -5, -2, 2 ),  & ( -3, 9, -4, -7, 6 ),  & ( -3, 9, -4, -5, 4 ),  & ( -3, 9, -4, -3, 2 ), \\
( -3, 9, -4, -1, 0 ),  & ( -3, 9, -3, -8, 6 ),  & ( -3, 9, -3, -6, 4 ),  & ( -3, 9, -3, -4, 2 ), \\
( -3, 9, -3, -2, 0 ),  & ( -3, 9, -2, -9, 6 ),  & ( -3, 9, -2, -7, 4 ),  & ( -3, 9, -2, -5, 2 ), \\
( -3, 9, -2, -3, 0 ),  & ( -3, 9, -1, -10, 6 ),  & ( -3, 9, -1, -8, 4 ),  & ( -3, 9, -1, -6, 2 ), \\
\end{array}}
$$

$$
{\tiny
\begin{array}{llll}
( -3, 9, -1, -4, 0 ),  & ( -3, 9, 0, -11, 6 ),  & ( -3, 9, 0, -9, 4 ),  & ( -3, 9, 0, -7, 2 ), \\
( -3, 9, 0, -5, 0 ),  & ( -3, 9, 1, -12, 6 ),  & ( -3, 9, 1, -10, 4 ),  & ( -3, 9, 1, -8, 2 ), \\
( -3, 9, 1, -6, 0 ),  & ( -3, 9, 2, -13, 6 ),  & ( -3, 9, 2, -11, 4 ),  & ( -3, 9, 2, -9, 2 ), \\
( -3, 9, 2, -7, 0 ),  & ( -3, 9, 3, -14, 6 ),  & ( -3, 9, 3, -10, 2 ),  & ( -3, 9, 3, -8, 0 ), \\
( -3, 9, 4, -15, 6 ),  & ( -3, 12, -7, -6, 5 ),  & ( -3, 12, -7, -4, 3 ),  & ( -3, 12, -6, -7, 5 ), \\
( -3, 12, -6, -5, 3 ),  & ( -3, 12, -5, -8, 5 ),  & ( -3, 12, -5, -6, 3 ),  & ( -3, 12, -4, -9, 5 ), \\
( -3, 12, -4, -7, 3 ),  & ( -3, 12, -3, -10, 5 ),  & ( -3, 12, -3, -8, 3 ),  & ( -3, 12, -2, -11, 5 ), \\
( -3, 12, -2, -9, 3 ),  & ( -3, 12, -1, -12, 5 ),  & ( -3, 12, -1, -10, 3 ),  & ( -3, 12, 0, -13, 5 ), \\
( -3, 12, 0, -11, 3 ),  & ( -3, 12, 1, -12, 3 ),  & ( -2, 0, 0, 4, -1 ),  & ( -2, 0, 1, 3, -1 ), \\
( -2, 0, 2, 2, -1 ),  & ( -2, 0, 2, 6, -5 ),  & ( -2, 0, 3, 1, -1 ),  & ( -2, 0, 3, 5, -5 ), \\
( -2, 0, 4, 0, -1 ),  & ( -2, 0, 4, 4, -5 ),  & ( -2, 0, 5, -1, -1 ),  & ( -2, 0, 5, 3, -5 ), \\
( -2, 0, 6, -2, -1 ),  & ( -2, 0, 6, 2, -5 ),  & ( -2, 0, 7, -3, -1 ),  & ( -2, 0, 7, 1, -5 ), \\
( -2, 0, 8, -4, -1 ),  & ( -2, 0, 8, 0, -5 ),  & ( -2, 3, -4, 2, 2 ),  & ( -2, 3, -3, 1, 2 ), \\
( -2, 3, -2, 0, 2 ),  & ( -2, 3, -2, 4, -2 ),  & ( -2, 3, -1, -1, 2 ),  & ( -2, 3, -1, 1, 0 ), \\
( -2, 3, -1, 3, -2 ),  & ( -2, 3, 0, -2, 2 ),  & ( -2, 3, 0, 0, 0 ),  & ( -2, 3, 0, 2, -2 ), \\
( -2, 3, 0, 4, -4 ),  & ( -2, 3, 1, -3, 2 ),  & ( -2, 3, 1, -1, 0 ),  & ( -2, 3, 1, 1, -2 ), \\
( -2, 3, 1, 3, -4 ),  & ( -2, 3, 2, -4, 2 ),  & ( -2, 3, 2, -2, 0 ),  & ( -2, 3, 2, 0, -2 ), \\
( -2, 3, 2, 2, -4 ),  & ( -2, 3, 3, -5, 2 ),  & ( -2, 3, 3, -3, 0 ),  & ( -2, 3, 3, -1, -2 ), \\
( -2, 3, 3, 1, -4 ),  & ( -2, 3, 4, -6, 2 ),  & ( -2, 3, 4, -4, 0 ),  & ( -2, 3, 4, -2, -2 ), \\
( -2, 3, 4, 0, -4 ),  & ( -2, 3, 5, -7, 2 ),  & ( -2, 3, 5, -3, -2 ),  & ( -2, 3, 5, -1, -4 ), \\
( -2, 3, 6, -8, 2 ),  & ( -2, 3, 6, -4, -2 ),  & ( -2, 3, 6, -2, -4 ),  & ( -2, 3, 7, -9, 2 ), \\
( -2, 3, 7, -5, -2 ),  & ( -2, 3, 8, -10, 2 ),  & ( -2, 3, 8, -6, -2 ),  & ( -2, 6, -7, -1, 5 ), \\
( -2, 6, -6, -2, 5 ),  & ( -2, 6, -6, 0, 3 ),  & ( -2, 6, -5, -3, 5 ),  & ( -2, 6, -5, -1, 3 ), \\
( -2, 6, -5, 1, 1 ),  & ( -2, 6, -4, -4, 5 ),  & ( -2, 6, -4, -2, 3 ),  & ( -2, 6, -4, 0, 1 ), \\
( -2, 6, -4, 2, -1 ),  & ( -2, 6, -3, -5, 5 ),  & ( -2, 6, -3, -3, 3 ),  & ( -2, 6, -3, -1, 1 ), \\
( -2, 6, -3, 1, -1 ),  & ( -2, 6, -2, -6, 5 ),  & ( -2, 6, -2, -4, 3 ),  & ( -2, 6, -2, -2, 1 ), \\
( -2, 6, -2, 0, -1 ),  & ( -2, 6, -1, -7, 5 ),  & ( -2, 6, -1, -5, 3 ),  & ( -2, 6, -1, -3, 1 ), \\
( -2, 6, -1, -1, -1 ),  & ( -2, 6, 0, -8, 5 ),  & ( -2, 6, 0, -6, 3 ),  & ( -2, 6, 0, -4, 1 ), \\
( -2, 6, 0, -2, -1 ),  & ( -2, 6, 1, -9, 5 ),  & ( -2, 6, 1, -7, 3 ),  & ( -2, 6, 1, -5, 1 ), \\
( -2, 6, 1, -3, -1 ),  & ( -2, 6, 2, -10, 5 ),  & ( -2, 6, 2, -8, 3 ),  & ( -2, 6, 2, -6, 1 ), \\
( -2, 6, 2, -4, -1 ),  & ( -2, 6, 3, -11, 5 ),  & ( -2, 6, 3, -9, 3 ),  & ( -2, 6, 3, -7, 1 ), \\
( -2, 6, 3, -5, -1 ),  & ( -2, 6, 4, -12, 5 ),  & ( -2, 6, 4, -10, 3 ),  & ( -2, 6, 4, -8, 1 ), \\
( -2, 6, 4, -6, -1 ),  & ( -2, 6, 5, -13, 5 ),  & ( -2, 6, 5, -9, 1 ),  & ( -2, 6, 5, -7, -1 ), \\
( -2, 6, 6, -10, 1 ),  & ( -2, 6, 6, -8, -1 ),  & ( -2, 9, -9, -5, 8 ),  & ( -2, 9, -9, -3, 6 ), \\
( -2, 9, -8, -6, 8 ),  & ( -2, 9, -8, -4, 6 ),  & ( -2, 9, -7, -7, 8 ),  & ( -2, 9, -7, -5, 6 ), \\
( -2, 9, -7, -3, 4 ),  & ( -2, 9, -7, -1, 2 ),  & ( -2, 9, -6, -8, 8 ),  & ( -2, 9, -6, -6, 6 ), \\
( -2, 9, -6, -4, 4 ),  & ( -2, 9, -6, -2, 2 ),  & ( -2, 9, -5, -9, 8 ),  & ( -2, 9, -5, -7, 6 ), \\
( -2, 9, -5, -5, 4 ),  & ( -2, 9, -5, -3, 2 ),  & ( -2, 9, -4, -10, 8 ),  & ( -2, 9, -4, -8, 6 ), \\
( -2, 9, -4, -6, 4 ),  & ( -2, 9, -4, -4, 2 ),  & ( -2, 9, -3, -11, 8 ),  & ( -2, 9, -3, -9, 6 ), \\
( -2, 9, -3, -7, 4 ),  & ( -2, 9, -3, -5, 2 ),  & ( -2, 9, -2, -12, 8 ),  & ( -2, 9, -2, -10, 6 ), \\
( -2, 9, -2, -8, 4 ),  & ( -2, 9, -2, -6, 2 ),  & ( -2, 9, -1, -13, 8 ),  & ( -2, 9, -1, -11, 6 ), \\
( -2, 9, -1, -9, 4 ),  & ( -2, 9, -1, -7, 2 ),  & ( -2, 9, 0, -14, 8 ),  & ( -2, 9, 0, -12, 6 ), \\
( -2, 9, 0, -10, 4 ),  & ( -2, 9, 0, -8, 2 ),  & ( -2, 9, 1, -15, 8 ),  & ( -2, 9, 1, -13, 6 ), \\
( -2, 9, 1, -11, 4 ),  & ( -2, 9, 1, -9, 2 ),  & ( -2, 9, 2, -14, 6 ),  & ( -2, 9, 2, -12, 4 ), \\
( -2, 9, 2, -10, 2 ),  & ( -2, 9, 3, -13, 4 ),  & ( -2, 9, 3, -11, 2 ),  & ( -2, 9, 4, -12, 2 ), \\
( -2, 12, -9, -7, 7 ),  & ( -2, 12, -9, -5, 5 ),  & ( -2, 12, -8, -8, 7 ),  & ( -2, 12, -8, -6, 5 ), \\
( -2, 12, -7, -9, 7 ),  & ( -2, 12, -7, -7, 5 ),  & ( -2, 12, -6, -10, 7 ),  & ( -2, 12, -6, -8, 5 ), \\
( -2, 12, -5, -11, 7 ),  & ( -2, 12, -5, -9, 5 ),  & ( -2, 12, -4, -12, 7 ),  & ( -2, 12, -4, -10, 5 ), \\
( -2, 12, -3, -13, 7 ),  & ( -2, 12, -3, -11, 5 ),  & ( -2, 12, -2, -12, 5 ),  & ( -2, 12, -1, -13, 5 ), \\
( -2, 12, 0, -14, 5 ),  & ( -1, -3, 0, 7, -2 ),  & ( -1, -3, 1, 6, -2 ),  & ( -1, -3, 2, 5, -2 ), \\
( -1, -3, 2, 9, -6 ),  & ( -1, -3, 3, 4, -2 ),  & ( -1, -3, 3, 8, -6 ),  & ( -1, -3, 4, 3, -2 ), \\
( -1, -3, 4, 7, -6 ),  & ( -1, -3, 5, 2, -2 ),  & ( -1, -3, 5, 6, -6 ),  & ( -1, -3, 6, 1, -2 ), \\
( -1, -3, 6, 5, -6 ),  & ( -1, -3, 7, 0, -2 ),  & ( -1, -3, 7, 4, -6 ),  & ( -1, -3, 8, -1, -2 ), \\
( -1, -3, 8, 3, -6 ),  & ( -1, -3, 9, -2, -2 ),  & ( -1, -3, 9, 2, -6 ),  & ( -1, -3, 10, -3, -2 ), \\
( -1, -3, 10, 1, -6 ),  & ( -1, -3, 11, 0, -6 ),  & ( -1, 0, -4, 5, 1 ),  & ( -1, 0, -3, 4, 1 ), \\
( -1, 0, -2, 3, 1 ),  & ( -1, 0, -2, 7, -3 ),  & ( -1, 0, -1, 2, 1 ),  & ( -1, 0, -1, 4, -1 ), \\
( -1, 0, -1, 6, -3 ),  & ( -1, 0, 0, 1, 1 ),  & ( -1, 0, 0, 3, -1 ),  & ( -1, 0, 0, 5, -3 ), \\
( -1, 0, 0, 7, -5 ),  & ( -1, 0, 1, 0, 1 ),  & ( -1, 0, 1, 2, -1 ),  & ( -1, 0, 1, 4, -3 ), \\
( -1, 0, 1, 6, -5 ),  & ( -1, 0, 2, -1, 1 ),  & ( -1, 0, 2, 1, -1 ),  & ( -1, 0, 2, 3, -3 ), \\
( -1, 0, 2, 5, -5 ),  & ( -1, 0, 3, -2, 1 ),  & ( -1, 0, 3, 0, -1 ),  & ( -1, 0, 3, 2, -3 ), \\
( -1, 0, 3, 4, -5 ),  & ( -1, 0, 4, -3, 1 ),  & ( -1, 0, 4, -1, -1 ),  & ( -1, 0, 4, 1, -3 ), \\
( -1, 0, 4, 3, -5 ),  & ( -1, 0, 5, -4, 1 ),  & ( -1, 0, 5, -2, -1 ),  & ( -1, 0, 5, 0, -3 ), \\
( -1, 0, 5, 2, -5 ),  & ( -1, 0, 6, -5, 1 ),  & ( -1, 0, 6, -3, -1 ),  & ( -1, 0, 6, -1, -3 ), \\
( -1, 0, 6, 1, -5 ),  & ( -1, 0, 7, -6, 1 ),  & ( -1, 0, 7, -2, -3 ),  & ( -1, 0, 7, 0, -5 ), \\
( -1, 0, 8, -7, 1 ),  & ( -1, 0, 8, -3, -3 ),  & ( -1, 0, 8, -1, -5 ),  & ( -1, 0, 9, -4, -3 ), \\
( -1, 3, -7, 2, 4 ),  & ( -1, 3, -6, 1, 4 ),  & ( -1, 3, -6, 3, 2 ),  & ( -1, 3, -5, 0, 4 ), \\
( -1, 3, -5, 2, 2 ),  & ( -1, 3, -5, 4, 0 ),  & ( -1, 3, -4, -1, 4 ),  & ( -1, 3, -4, 1, 2 ), \\
( -1, 3, -4, 3, 0 ),  & ( -1, 3, -4, 5, -2 ),  & ( -1, 3, -3, -2, 4 ),  & ( -1, 3, -3, 0, 2 ), \\
( -1, 3, -3, 2, 0 ),  & ( -1, 3, -3, 4, -2 ),  & ( -1, 3, -2, -3, 4 ),  & ( -1, 3, -2, -1, 2 ), \\
( -1, 3, -2, 1, 0 ),  & ( -1, 3, -2, 3, -2 ),  & ( -1, 3, -1, -4, 4 ),  & ( -1, 3, -1, -2, 2 ), \\
( -1, 3, -1, 0, 0 ),  & ( -1, 3, -1, 2, -2 ),  & ( -1, 3, 0, -5, 4 ),  & ( -1, 3, 0, -3, 2 ), \\
( -1, 3, 0, -1, 0 ),  & ( -1, 3, 0, 1, -2 ),  & ( -1, 3, 1, -6, 4 ),  & ( -1, 3, 1, -4, 2 ), \\
( -1, 3, 1, -2, 0 ),  & ( -1, 3, 1, 0, -2 ),  & ( -1, 3, 2, -7, 4 ),  & ( -1, 3, 2, -5, 2 ), \\
( -1, 3, 2, -3, 0 ),  & ( -1, 3, 2, -1, -2 ),  & ( -1, 3, 3, -8, 4 ),  & ( -1, 3, 3, -6, 2 ), \\
( -1, 3, 3, -4, 0 ),  & ( -1, 3, 3, -2, -2 ),  & ( -1, 3, 4, -9, 4 ),  & ( -1, 3, 4, -7, 2 ), \\
( -1, 3, 4, -5, 0 ),  & ( -1, 3, 4, -3, -2 ),  & ( -1, 3, 5, -10, 4 ),  & ( -1, 3, 5, -8, 2 ), \\
( -1, 3, 5, -6, 0 ),  & ( -1, 3, 5, -4, -2 ),  & ( -1, 3, 6, -9, 2 ),  & ( -1, 3, 6, -7, 0 ), \\
( -1, 3, 6, -5, -2 ),  & ( -1, 3, 7, -8, 0 ),  & ( -1, 3, 7, -6, -2 ),  & ( -1, 3, 8, -7, -2 ), \\
( -1, 6, -9, -2, 7 ),  & ( -1, 6, -9, 0, 5 ),  & ( -1, 6, -8, -3, 7 ),  & ( -1, 6, -8, -1, 5 ), \\
( -1, 6, -7, -4, 7 ),  & ( -1, 6, -7, -2, 5 ),  & ( -1, 6, -7, 0, 3 ),  & ( -1, 6, -7, 2, 1 ), \\
\end{array}}
$$

$$
{\tiny
\begin{array}{llll}
( -1, 6, -6, -5, 7 ),  & ( -1, 6, -6, -3, 5 ),  & ( -1, 6, -6, -1, 3 ),  & ( -1, 6, -6, 1, 1 ), \\
( -1, 6, -5, -6, 7 ),  & ( -1, 6, -5, -4, 5 ),  & ( -1, 6, -5, -2, 3 ),  & ( -1, 6, -5, 0, 1 ), \\
( -1, 6, -4, -7, 7 ),  & ( -1, 6, -4, -5, 5 ),  & ( -1, 6, -4, -3, 3 ),  & ( -1, 6, -4, -1, 1 ), \\
( -1, 6, -3, -8, 7 ),  & ( -1, 6, -3, -6, 5 ),  & ( -1, 6, -3, -4, 3 ),  & ( -1, 6, -3, -2, 1 ), \\
( -1, 6, -2, -9, 7 ),  & ( -1, 6, -2, -7, 5 ),  & ( -1, 6, -2, -5, 3 ),  & ( -1, 6, -2, -3, 1 ), \\
( -1, 6, -1, -10, 7 ),  & ( -1, 6, -1, -8, 5 ),  & ( -1, 6, -1, -6, 3 ),  & ( -1, 6, -1, -4, 1 ), \\
( -1, 6, 0, -11, 7 ),  & ( -1, 6, 0, -9, 5 ),  & ( -1, 6, 0, -7, 3 ),  & ( -1, 6, 0, -5, 1 ), \\
( -1, 6, 1, -12, 7 ),  & ( -1, 6, 1, -10, 5 ),  & ( -1, 6, 1, -8, 3 ),  & ( -1, 6, 1, -6, 1 ), \\
( -1, 6, 2, -11, 5 ),  & ( -1, 6, 2, -9, 3 ),  & ( -1, 6, 2, -7, 1 ),  & ( -1, 6, 3, -12, 5 ), \\
( -1, 6, 3, -10, 3 ),  & ( -1, 6, 3, -8, 1 ),  & ( -1, 6, 4, -9, 1 ),  & ( -1, 6, 5, -10, 1 ), \\
( -1, 9, -10, -7, 10 ),  & ( -1, 9, -10, -5, 8 ),  & ( -1, 9, -9, -8, 10 ),  & ( -1, 9, -9, -6, 8 ), \\
( -1, 9, -9, -4, 6 ),  & ( -1, 9, -9, -2, 4 ),  & ( -1, 9, -8, -9, 10 ),  & ( -1, 9, -8, -7, 8 ), \\
( -1, 9, -8, -5, 6 ),  & ( -1, 9, -8, -3, 4 ),  & ( -1, 9, -7, -10, 10 ),  & ( -1, 9, -7, -8, 8 ), \\
( -1, 9, -7, -6, 6 ),  & ( -1, 9, -7, -4, 4 ),  & ( -1, 9, -6, -11, 10 ),  & ( -1, 9, -6, -9, 8 ), \\
( -1, 9, -6, -7, 6 ),  & ( -1, 9, -6, -5, 4 ),  & ( -1, 9, -5, -12, 10 ),  & ( -1, 9, -5, -10, 8 ), \\
( -1, 9, -5, -8, 6 ),  & ( -1, 9, -5, -6, 4 ),  & ( -1, 9, -4, -13, 10 ),  & ( -1, 9, -4, -11, 8 ), \\
( -1, 9, -4, -9, 6 ),  & ( -1, 9, -4, -7, 4 ),  & ( -1, 9, -3, -14, 10 ),  & ( -1, 9, -3, -12, 8 ), \\
( -1, 9, -3, -10, 6 ),  & ( -1, 9, -3, -8, 4 ),  & ( -1, 9, -2, -13, 8 ),  & ( -1, 9, -2, -9, 4 ), \\
( -1, 9, -1, -14, 8 ),  & ( -1, 9, -1, -10, 4 ),  & ( -1, 9, 0, -11, 4 ),  & ( -1, 9, 1, -12, 4 ), \\
( -1, 12, -10, -7, 7 ),  & ( -1, 12, -9, -8, 7 ),  & ( -1, 12, -8, -9, 7 ),  & ( -1, 12, -7, -10, 7 ), \\
( -1, 12, -6, -11, 7 ),  & ( -1, 12, -5, -12, 7 ),  & ( -1, 12, -4, -13, 7 ),  & ( -1, 12, -3, -14, 7 ), \\
( 0, -9, 5, 11, -6 ),  & ( 0, -9, 6, 10, -6 ),  & ( 0, -9, 7, 9, -6 ),  & ( 0, -9, 8, 8, -6 ), \\
( 0, -9, 9, 7, -6 ),  & ( 0, -9, 10, 6, -6 ),  & ( 0, -9, 11, 5, -6 ),  & ( 0, -9, 12, 4, -6 ), \\
( 0, -6, 0, 10, -3 ),  & ( 0, -6, 1, 9, -3 ),  & ( 0, -6, 2, 8, -3 ),  & ( 0, -6, 2, 12, -7 ), \\
( 0, -6, 3, 7, -3 ),  & ( 0, -6, 3, 11, -7 ),  & ( 0, -6, 4, 6, -3 ),  & ( 0, -6, 4, 10, -7 ), \\
( 0, -6, 5, 5, -3 ),  & ( 0, -6, 5, 9, -7 ),  & ( 0, -6, 6, 4, -3 ),  & ( 0, -6, 6, 8, -7 ), \\
( 0, -6, 7, 3, -3 ),  & ( 0, -6, 7, 7, -7 ),  & ( 0, -6, 8, 2, -3 ),  & ( 0, -6, 8, 6, -7 ), \\
( 0, -6, 9, 1, -3 ),  & ( 0, -6, 9, 5, -7 ),  & ( 0, -6, 10, 0, -3 ),  & ( 0, -6, 10, 4, -7 ), \\
( 0, -6, 11, 3, -7 ),  & ( 0, -6, 12, 2, -7 ),  & ( 0, -3, -4, 8, 0 ),  & ( 0, -3, -3, 7, 0 ), \\
( 0, -3, -2, 6, 0 ),  & ( 0, -3, -2, 10, -4 ),  & ( 0, -3, -1, 5, 0 ),  & ( 0, -3, -1, 7, -2 ), \\
( 0, -3, -1, 9, -4 ),  & ( 0, -3, 0, 4, 0 ),  & ( 0, -3, 0, 6, -2 ),  & ( 0, -3, 0, 8, -4 ), \\
( 0, -3, 0, 10, -6 ),  & ( 0, -3, 1, 3, 0 ),  & ( 0, -3, 1, 5, -2 ),  & ( 0, -3, 1, 7, -4 ), \\
( 0, -3, 1, 9, -6 ),  & ( 0, -3, 2, 2, 0 ),  & ( 0, -3, 2, 4, -2 ),  & ( 0, -3, 2, 6, -4 ), \\
( 0, -3, 2, 8, -6 ),  & ( 0, -3, 3, 1, 0 ),  & ( 0, -3, 3, 3, -2 ),  & ( 0, -3, 3, 5, -4 ), \\
( 0, -3, 3, 7, -6 ),  & ( 0, -3, 4, 0, 0 ),  & ( 0, -3, 4, 2, -2 ),  & ( 0, -3, 4, 4, -4 ), \\
( 0, -3, 4, 6, -6 ),  & ( 0, -3, 5, -1, 0 ),  & ( 0, -3, 5, 1, -2 ),  & ( 0, -3, 5, 3, -4 ), \\
( 0, -3, 5, 5, -6 ),  & ( 0, -3, 6, -2, 0 ),  & ( 0, -3, 6, 0, -2 ),  & ( 0, -3, 6, 2, -4 ), \\
( 0, -3, 6, 4, -6 ),  & ( 0, -3, 7, -3, 0 ),  & ( 0, -3, 7, -1, -2 ),  & ( 0, -3, 7, 1, -4 ), \\
( 0, -3, 7, 3, -6 ),  & ( 0, -3, 8, -4, 0 ),  & ( 0, -3, 8, -2, -2 ),  & ( 0, -3, 8, 0, -4 ), \\
( 0, -3, 8, 2, -6 ),  & ( 0, -3, 9, -1, -4 ),  & ( 0, -3, 9, 1, -6 ),  & ( 0, -3, 10, -2, -4 ), \\
( 0, -3, 10, 0, -6 ),  & ( 0, 0, -7, 5, 3 ),  & ( 0, 0, -6, 4, 3 ),  & ( 0, 0, -6, 6, 1 ), \\
( 0, 0, -5, 3, 3 ),  & ( 0, 0, -5, 5, 1 ),  & ( 0, 0, -5, 7, -1 ),  & ( 0, 0, -4, 2, 3 ), \\
( 0, 0, -4, 4, 1 ),  & ( 0, 0, -4, 6, -1 ),  & ( 0, 0, -4, 8, -3 ),  & ( 0, 0, -3, 1, 3 ), \\
( 0, 0, -3, 3, 1 ),  & ( 0, 0, -3, 5, -1 ),  & ( 0, 0, -3, 7, -3 ),  & ( 0, 0, -2, 0, 3 ), \\
( 0, 0, -2, 2, 1 ),  & ( 0, 0, -2, 4, -1 ),  & ( 0, 0, -2, 6, -3 ),  & ( 0, 0, -1, -1, 3 ), \\
( 0, 0, -1, 1, 1 ),  & ( 0, 0, -1, 3, -1 ),  & ( 0, 0, -1, 5, -3 ),  & ( 0, 0, 0, -2, 3 ), \\
( 0, 0, 0, 0, 1 ),  & ( 0, 0, 0, 2, -1 ),  & ( 0, 0, 0, 4, -3 ),  & ( 0, 0, 1, -3, 3 ), \\
( 0, 0, 1, -1, 1 ),  & ( 0, 0, 1, 1, -1 ),  & ( 0, 0, 1, 3, -3 ),  & ( 0, 0, 2, -4, 3 ), \\
( 0, 0, 2, -2, 1 ),  & ( 0, 0, 2, 0, -1 ),  & ( 0, 0, 2, 2, -3 ),  & ( 0, 0, 3, -5, 3 ), \\
( 0, 0, 3, -3, 1 ),  & ( 0, 0, 3, -1, -1 ),  & ( 0, 0, 3, 1, -3 ),  & ( 0, 0, 4, -6, 3 ), \\
( 0, 0, 4, -4, 1 ),  & ( 0, 0, 4, -2, -1 ),  & ( 0, 0, 4, 0, -3 ),  & ( 0, 0, 5, -7, 3 ), \\
( 0, 0, 5, -5, 1 ),  & ( 0, 0, 5, -3, -1 ),  & ( 0, 0, 5, -1, -3 ),  & ( 0, 0, 6, -6, 1 ), \\
( 0, 0, 6, -4, -1 ),  & ( 0, 0, 6, -2, -3 ),  & ( 0, 0, 7, -5, -1 ),  & ( 0, 0, 7, -3, -3 ), \\
( 0, 0, 8, -4, -3 ),  & ( 0, 3, -9, 1, 6 ),  & ( 0, 3, -9, 3, 4 ),  & ( 0, 3, -8, 0, 6 ), \\
( 0, 3, -8, 2, 4 ),  & ( 0, 3, -7, -1, 6 ),  & ( 0, 3, -7, 1, 4 ),  & ( 0, 3, -7, 3, 2 ), \\
( 0, 3, -7, 5, 0 ),  & ( 0, 3, -6, -2, 6 ),  & ( 0, 3, -6, 0, 4 ),  & ( 0, 3, -6, 2, 2 ), \\
( 0, 3, -6, 4, 0 ),  & ( 0, 3, -5, -3, 6 ),  & ( 0, 3, -5, -1, 4 ),  & ( 0, 3, -5, 1, 2 ), \\
( 0, 3, -5, 3, 0 ),  & ( 0, 3, -4, -4, 6 ),  & ( 0, 3, -4, -2, 4 ),  & ( 0, 3, -4, 0, 2 ), \\
( 0, 3, -4, 2, 0 ),  & ( 0, 3, -3, -5, 6 ),  & ( 0, 3, -3, -3, 4 ),  & ( 0, 3, -3, -1, 2 ), \\
( 0, 3, -3, 1, 0 ),  & ( 0, 3, -2, -6, 6 ),  & ( 0, 3, -2, -4, 4 ),  & ( 0, 3, -2, -2, 2 ), \\
( 0, 3, -2, 0, 0 ),  & ( 0, 3, -1, -7, 6 ),  & ( 0, 3, -1, -5, 4 ),  & ( 0, 3, -1, -3, 2 ), \\
( 0, 3, -1, -1, 0 ),  & ( 0, 3, 0, -8, 6 ),  & ( 0, 3, 0, -6, 4 ),  & ( 0, 3, 0, -4, 2 ), \\
( 0, 3, 0, -2, 0 ),  & ( 0, 3, 1, -9, 6 ),  & ( 0, 3, 1, -7, 4 ),  & ( 0, 3, 1, -5, 2 ), \\
( 0, 3, 1, -3, 0 ),  & ( 0, 3, 2, -8, 4 ),  & ( 0, 3, 2, -6, 2 ),  & ( 0, 3, 2, -4, 0 ), \\
( 0, 3, 3, -9, 4 ),  & ( 0, 3, 3, -7, 2 ),  & ( 0, 3, 3, -5, 0 ),  & ( 0, 3, 4, -6, 0 ), \\
( 0, 3, 5, -7, 0 ),  & ( 0, 6, -10, -2, 7 ),  & ( 0, 6, -9, -5, 9 ),  & ( 0, 6, -9, -3, 7 ), \\
( 0, 6, -9, 1, 3 ),  & ( 0, 6, -8, -6, 9 ),  & ( 0, 6, -8, -4, 7 ),  & ( 0, 6, -8, -2, 5 ), \\
( 0, 6, -8, 0, 3 ),  & ( 0, 6, -7, -7, 9 ),  & ( 0, 6, -7, -5, 7 ),  & ( 0, 6, -7, -3, 5 ), \\
( 0, 6, -7, -1, 3 ),  & ( 0, 6, -6, -8, 9 ),  & ( 0, 6, -6, -6, 7 ),  & ( 0, 6, -6, -4, 5 ), \\
( 0, 6, -6, -2, 3 ),  & ( 0, 6, -5, -9, 9 ),  & ( 0, 6, -5, -7, 7 ),  & ( 0, 6, -5, -5, 5 ), \\
( 0, 6, -5, -3, 3 ),  & ( 0, 6, -4, -10, 9 ),  & ( 0, 6, -4, -8, 7 ),  & ( 0, 6, -4, -6, 5 ), \\
( 0, 6, -4, -4, 3 ),  & ( 0, 6, -3, -11, 9 ),  & ( 0, 6, -3, -9, 7 ),  & ( 0, 6, -3, -7, 5 ), \\
( 0, 6, -3, -5, 3 ),  & ( 0, 6, -2, -10, 7 ),  & ( 0, 6, -2, -6, 3 ),  & ( 0, 6, -1, -11, 7 ), \\
( 0, 6, -1, -7, 3 ),  & ( 0, 6, 0, -8, 3 ),  & ( 0, 6, 1, -9, 3 ),  & ( 0, 9, -10, -4, 6 ), \\
( 0, 9, -9, -5, 6 ),  & ( 0, 9, -8, -6, 6 ),  & ( 0, 9, -7, -7, 6 ),  & ( 0, 9, -6, -8, 6 ), \\
( 0, 9, -5, -9, 6 ),  & ( 0, 9, -4, -10, 6 ),  & ( 0, 9, -3, -11, 6 ),  & ( 1, -9, 0, 13, -4 ), \\
( 1, -9, 1, 12, -4 ),  & ( 1, -9, 2, 11, -4 ),  & ( 1, -9, 3, 10, -4 ),  & ( 1, -9, 4, 9, -4 ), \\
( 1, -9, 5, 8, -4 ),  & ( 1, -9, 5, 10, -6 ),  & ( 1, -9, 6, 7, -4 ),  & ( 1, -9, 6, 9, -6 ), \\
( 1, -9, 7, 6, -4 ),  & ( 1, -9, 7, 8, -6 ),  & ( 1, -9, 8, 5, -4 ),  & ( 1, -9, 8, 7, -6 ), \\
( 1, -9, 9, 4, -4 ),  & ( 1, -9, 9, 6, -6 ),  & ( 1, -9, 10, 3, -4 ),  & ( 1, -9, 10, 5, -6 ), \\
( 1, -6, -4, 11, -1 ),  & ( 1, -6, -3, 10, -1 ),  & ( 1, -6, -2, 9, -1 ),  & ( 1, -6, -2, 13, -5 ), \\
( 1, -6, -1, 8, -1 ),  & ( 1, -6, -1, 10, -3 ),  & ( 1, -6, -1, 12, -5 ),  & ( 1, -6, 0, 7, -1 ), \\
( 1, -6, 0, 9, -3 ),  & ( 1, -6, 0, 11, -5 ),  & ( 1, -6, 0, 13, -7 ),  & ( 1, -6, 1, 6, -1 ), \\
( 1, -6, 1, 8, -3 ),  & ( 1, -6, 1, 10, -5 ),  & ( 1, -6, 1, 12, -7 ),  & ( 1, -6, 2, 5, -1 ), \\
( 1, -6, 2, 7, -3 ),  & ( 1, -6, 2, 9, -5 ),  & ( 1, -6, 2, 11, -7 ),  & ( 1, -6, 3, 4, -1 ), \\
( 1, -6, 3, 6, -3 ),  & ( 1, -6, 3, 8, -5 ),  & ( 1, -6, 3, 10, -7 ),  & ( 1, -6, 4, 3, -1 ), \\
\end{array}}
$$

$$
{\tiny
\begin{array}{llll}
( 1, -6, 4, 5, -3 ),  & ( 1, -6, 4, 7, -5 ),  & ( 1, -6, 4, 9, -7 ),  & ( 1, -6, 5, 2, -1 ), \\
( 1, -6, 5, 4, -3 ),  & ( 1, -6, 5, 6, -5 ),  & ( 1, -6, 5, 8, -7 ),  & ( 1, -6, 6, 1, -1 ), \\
( 1, -6, 6, 3, -3 ),  & ( 1, -6, 6, 5, -5 ),  & ( 1, -6, 6, 7, -7 ),  & ( 1, -6, 7, 0, -1 ), \\
( 1, -6, 7, 2, -3 ),  & ( 1, -6, 7, 4, -5 ),  & ( 1, -6, 7, 6, -7 ),  & ( 1, -6, 8, -1, -1 ), \\
( 1, -6, 8, 1, -3 ),  & ( 1, -6, 8, 3, -5 ),  & ( 1, -6, 8, 5, -7 ),  & ( 1, -6, 9, 2, -5 ), \\
( 1, -6, 9, 4, -7 ),  & ( 1, -6, 10, 1, -5 ),  & ( 1, -6, 10, 3, -7 ),  & ( 1, -3, -7, 8, 2 ), \\
( 1, -3, -6, 7, 2 ),  & ( 1, -3, -6, 9, 0 ),  & ( 1, -3, -5, 6, 2 ),  & ( 1, -3, -5, 8, 0 ), \\
( 1, -3, -5, 10, -2 ),  & ( 1, -3, -4, 5, 2 ),  & ( 1, -3, -4, 7, 0 ),  & ( 1, -3, -4, 9, -2 ), \\
( 1, -3, -4, 11, -4 ),  & ( 1, -3, -3, 4, 2 ),  & ( 1, -3, -3, 6, 0 ),  & ( 1, -3, -3, 8, -2 ), \\
( 1, -3, -3, 10, -4 ),  & ( 1, -3, -2, 3, 2 ),  & ( 1, -3, -2, 5, 0 ),  & ( 1, -3, -2, 7, -2 ), \\
( 1, -3, -2, 9, -4 ),  & ( 1, -3, -1, 2, 2 ),  & ( 1, -3, -1, 4, 0 ),  & ( 1, -3, -1, 6, -2 ), \\
( 1, -3, -1, 8, -4 ),  & ( 1, -3, 0, 1, 2 ),  & ( 1, -3, 0, 3, 0 ),  & ( 1, -3, 0, 5, -2 ), \\
( 1, -3, 0, 7, -4 ),  & ( 1, -3, 1, 0, 2 ),  & ( 1, -3, 1, 2, 0 ),  & ( 1, -3, 1, 4, -2 ), \\
( 1, -3, 1, 6, -4 ),  & ( 1, -3, 2, -1, 2 ),  & ( 1, -3, 2, 1, 0 ),  & ( 1, -3, 2, 3, -2 ), \\
( 1, -3, 2, 5, -4 ),  & ( 1, -3, 3, -2, 2 ),  & ( 1, -3, 3, 0, 0 ),  & ( 1, -3, 3, 2, -2 ), \\
( 1, -3, 3, 4, -4 ),  & ( 1, -3, 4, -3, 2 ),  & ( 1, -3, 4, -1, 0 ),  & ( 1, -3, 4, 1, -2 ), \\
( 1, -3, 4, 3, -4 ),  & ( 1, -3, 5, -4, 2 ),  & ( 1, -3, 5, -2, 0 ),  & ( 1, -3, 5, 0, -2 ), \\
( 1, -3, 5, 2, -4 ),  & ( 1, -3, 6, -3, 0 ),  & ( 1, -3, 6, -1, -2 ),  & ( 1, -3, 6, 1, -4 ), \\
( 1, -3, 7, -2, -2 ),  & ( 1, -3, 7, 0, -4 ),  & ( 1, -3, 8, -1, -4 ),  & ( 1, 0, -8, 3, 5 ), \\
( 1, 0, -8, 5, 3 ),  & ( 1, 0, -7, 2, 5 ),  & ( 1, 0, -7, 4, 3 ),  & ( 1, 0, -7, 8, -1 ), \\
( 1, 0, -6, 1, 5 ),  & ( 1, 0, -6, 3, 3 ),  & ( 1, 0, -6, 5, 1 ),  & ( 1, 0, -6, 7, -1 ), \\
( 1, 0, -5, 0, 5 ),  & ( 1, 0, -5, 2, 3 ),  & ( 1, 0, -5, 4, 1 ),  & ( 1, 0, -5, 6, -1 ), \\
( 1, 0, -4, -1, 5 ),  & ( 1, 0, -4, 1, 3 ),  & ( 1, 0, -4, 3, 1 ),  & ( 1, 0, -4, 5, -1 ), \\
( 1, 0, -3, -2, 5 ),  & ( 1, 0, -3, 0, 3 ),  & ( 1, 0, -3, 2, 1 ),  & ( 1, 0, -3, 4, -1 ), \\
( 1, 0, -2, -3, 5 ),  & ( 1, 0, -2, -1, 3 ),  & ( 1, 0, -2, 1, 1 ),  & ( 1, 0, -2, 3, -1 ), \\
( 1, 0, -1, -4, 5 ),  & ( 1, 0, -1, -2, 3 ),  & ( 1, 0, -1, 0, 1 ),  & ( 1, 0, -1, 2, -1 ), \\
( 1, 0, 0, -5, 5 ),  & ( 1, 0, 0, -3, 3 ),  & ( 1, 0, 0, -1, 1 ),  & ( 1, 0, 0, 1, -1 ), \\
( 1, 0, 1, -6, 5 ),  & ( 1, 0, 1, -4, 3 ),  & ( 1, 0, 1, -2, 1 ),  & ( 1, 0, 1, 0, -1 ), \\
( 1, 0, 2, -5, 3 ),  & ( 1, 0, 2, -3, 1 ),  & ( 1, 0, 2, -1, -1 ),  & ( 1, 0, 3, -6, 3 ), \\
( 1, 0, 3, -4, 1 ),  & ( 1, 0, 3, -2, -1 ),  & ( 1, 0, 4, -3, -1 ),  & ( 1, 0, 5, -4, -1 ), \\
( 1, 3, -10, 1, 6 ),  & ( 1, 3, -9, 0, 6 ),  & ( 1, 3, -9, 4, 2 ),  & ( 1, 3, -8, -1, 6 ), \\
( 1, 3, -8, 3, 2 ),  & ( 1, 3, -7, -2, 6 ),  & ( 1, 3, -7, 2, 2 ),  & ( 1, 3, -6, -3, 6 ), \\
( 1, 3, -6, 1, 2 ),  & ( 1, 3, -5, -4, 6 ),  & ( 1, 3, -5, 0, 2 ),  & ( 1, 3, -4, -5, 6 ), \\
( 1, 3, -4, -1, 2 ),  & ( 1, 3, -3, -6, 6 ),  & ( 1, 3, -3, -2, 2 ),  & ( 1, 3, -2, -7, 6 ), \\
( 1, 3, -2, -3, 2 ),  & ( 1, 3, -1, -8, 6 ),  & ( 1, 3, -1, -4, 2 ),  & ( 1, 3, 0, -5, 2 ), \\
( 1, 3, 1, -6, 2 ),  & ( 1, 6, -9, -2, 5 ),  & ( 1, 6, -8, -3, 5 ),  & ( 1, 6, -7, -4, 5 ), \\
( 1, 6, -6, -5, 5 ),  & ( 1, 6, -5, -6, 5 ),  & ( 1, 6, -4, -7, 5 ),  & ( 1, 6, -3, -8, 5 ), \\
( 2, -12, 3, 13, -5 ),  & ( 2, -12, 4, 12, -5 ),  & ( 2, -12, 5, 11, -5 ),  & ( 2, -12, 6, 10, -5 ), \\
( 2, -12, 7, 9, -5 ),  & ( 2, -12, 8, 8, -5 ),  & ( 2, -12, 9, 7, -5 ),  & ( 2, -9, -4, 14, -2 ), \\
( 2, -9, -3, 13, -2 ),  & ( 2, -9, -2, 12, -2 ),  & ( 2, -9, -1, 11, -2 ),  & ( 2, -9, -1, 13, -4 ), \\
( 2, -9, 0, 10, -2 ),  & ( 2, -9, 0, 12, -4 ),  & ( 2, -9, 0, 14, -6 ),  & ( 2, -9, 1, 9, -2 ), \\
( 2, -9, 1, 11, -4 ),  & ( 2, -9, 1, 13, -6 ),  & ( 2, -9, 2, 8, -2 ),  & ( 2, -9, 2, 10, -4 ), \\
( 2, -9, 2, 12, -6 ),  & ( 2, -9, 2, 14, -8 ),  & ( 2, -9, 3, 7, -2 ),  & ( 2, -9, 3, 9, -4 ), \\
( 2, -9, 3, 11, -6 ),  & ( 2, -9, 3, 13, -8 ),  & ( 2, -9, 4, 6, -2 ),  & ( 2, -9, 4, 8, -4 ), \\
( 2, -9, 4, 10, -6 ),  & ( 2, -9, 4, 12, -8 ),  & ( 2, -9, 5, 5, -2 ),  & ( 2, -9, 5, 7, -4 ), \\
( 2, -9, 5, 9, -6 ),  & ( 2, -9, 5, 11, -8 ),  & ( 2, -9, 6, 4, -2 ),  & ( 2, -9, 6, 6, -4 ), \\
( 2, -9, 6, 8, -6 ),  & ( 2, -9, 6, 10, -8 ),  & ( 2, -9, 7, 3, -2 ),  & ( 2, -9, 7, 5, -4 ), \\
( 2, -9, 7, 7, -6 ),  & ( 2, -9, 7, 9, -8 ),  & ( 2, -9, 8, 2, -2 ),  & ( 2, -9, 8, 4, -4 ), \\
( 2, -9, 8, 6, -6 ),  & ( 2, -9, 8, 8, -8 ),  & ( 2, -9, 9, 5, -6 ),  & ( 2, -9, 9, 7, -8 ), \\
( 2, -6, -6, 10, 1 ),  & ( 2, -6, -5, 9, 1 ),  & ( 2, -6, -5, 11, -1 ),  & ( 2, -6, -4, 8, 1 ), \\
( 2, -6, -4, 10, -1 ),  & ( 2, -6, -4, 12, -3 ),  & ( 2, -6, -4, 14, -5 ),  & ( 2, -6, -3, 7, 1 ), \\
( 2, -6, -3, 9, -1 ),  & ( 2, -6, -3, 11, -3 ),  & ( 2, -6, -3, 13, -5 ),  & ( 2, -6, -2, 6, 1 ), \\
( 2, -6, -2, 8, -1 ),  & ( 2, -6, -2, 10, -3 ),  & ( 2, -6, -2, 12, -5 ),  & ( 2, -6, -1, 5, 1 ), \\
( 2, -6, -1, 7, -1 ),  & ( 2, -6, -1, 9, -3 ),  & ( 2, -6, -1, 11, -5 ),  & ( 2, -6, 0, 4, 1 ), \\
( 2, -6, 0, 6, -1 ),  & ( 2, -6, 0, 8, -3 ),  & ( 2, -6, 0, 10, -5 ),  & ( 2, -6, 1, 3, 1 ), \\
( 2, -6, 1, 5, -1 ),  & ( 2, -6, 1, 7, -3 ),  & ( 2, -6, 1, 9, -5 ),  & ( 2, -6, 2, 2, 1 ), \\
( 2, -6, 2, 4, -1 ),  & ( 2, -6, 2, 6, -3 ),  & ( 2, -6, 2, 8, -5 ),  & ( 2, -6, 3, 1, 1 ), \\
( 2, -6, 3, 3, -1 ),  & ( 2, -6, 3, 5, -3 ),  & ( 2, -6, 3, 7, -5 ),  & ( 2, -6, 4, 0, 1 ), \\
( 2, -6, 4, 2, -1 ),  & ( 2, -6, 4, 4, -3 ),  & ( 2, -6, 4, 6, -5 ),  & ( 2, -6, 5, -1, 1 ), \\
( 2, -6, 5, 1, -1 ),  & ( 2, -6, 5, 3, -3 ),  & ( 2, -6, 5, 5, -5 ),  & ( 2, -6, 6, 0, -1 ), \\
( 2, -6, 6, 2, -3 ),  & ( 2, -6, 6, 4, -5 ),  & ( 2, -6, 7, 1, -3 ),  & ( 2, -6, 7, 3, -5 ), \\
( 2, -6, 8, 2, -5 ),  & ( 2, -3, -7, 7, 2 ),  & ( 2, -3, -7, 11, -2 ),  & ( 2, -3, -6, 4, 4 ), \\
( 2, -3, -6, 6, 2 ),  & ( 2, -3, -6, 10, -2 ),  & ( 2, -3, -5, 3, 4 ),  & ( 2, -3, -5, 5, 2 ), \\
( 2, -3, -5, 9, -2 ),  & ( 2, -3, -4, 2, 4 ),  & ( 2, -3, -4, 4, 2 ),  & ( 2, -3, -4, 6, 0 ), \\
( 2, -3, -4, 8, -2 ),  & ( 2, -3, -3, 1, 4 ),  & ( 2, -3, -3, 3, 2 ),  & ( 2, -3, -3, 5, 0 ), \\
( 2, -3, -3, 7, -2 ),  & ( 2, -3, -2, 0, 4 ),  & ( 2, -3, -2, 2, 2 ),  & ( 2, -3, -2, 4, 0 ), \\
( 2, -3, -2, 6, -2 ),  & ( 2, -3, -1, -1, 4 ),  & ( 2, -3, -1, 1, 2 ),  & ( 2, -3, -1, 3, 0 ), \\
( 2, -3, -1, 5, -2 ),  & ( 2, -3, 0, -2, 4 ),  & ( 2, -3, 0, 0, 2 ),  & ( 2, -3, 0, 2, 0 ), \\
( 2, -3, 0, 4, -2 ),  & ( 2, -3, 1, -3, 4 ),  & ( 2, -3, 1, -1, 2 ),  & ( 2, -3, 1, 1, 0 ), \\
( 2, -3, 1, 3, -2 ),  & ( 2, -3, 2, -2, 2 ),  & ( 2, -3, 2, 0, 0 ),  & ( 2, -3, 2, 2, -2 ), \\
( 2, -3, 3, -3, 2 ),  & ( 2, -3, 3, -1, 0 ),  & ( 2, -3, 3, 1, -2 ),  & ( 2, -3, 4, 0, -2 ), \\
( 2, -3, 5, -1, -2 ),  & ( 2, 0, -8, 2, 5 ),  & ( 2, 0, -8, 6, 1 ),  & ( 2, 0, -7, 1, 5 ), \\
( 2, 0, -7, 5, 1 ),  & ( 2, 0, -6, 0, 5 ),  & ( 2, 0, -6, 4, 1 ),  & ( 2, 0, -5, -1, 5 ), \\
( 2, 0, -5, 3, 1 ),  & ( 2, 0, -4, -2, 5 ),  & ( 2, 0, -4, 2, 1 ),  & ( 2, 0, -3, -3, 5 ), \\
( 2, 0, -3, 1, 1 ),  & ( 2, 0, -2, -4, 5 ),  & ( 2, 0, -2, 0, 1 ),  & ( 2, 0, -1, -5, 5 ), \\
( 2, 0, -1, -1, 1 ),  & ( 2, 0, 0, -2, 1 ),  & ( 2, 0, 1, -3, 1 ),  & ( 3, -12, 0, 13, -3 ), \\
( 3, -12, 1, 12, -3 ),  & ( 3, -12, 2, 11, -3 ),  & ( 3, -12, 3, 10, -3 ),  & ( 3, -12, 3, 12, -5 ), \\
( 3, -12, 4, 9, -3 ),  & ( 3, -12, 4, 11, -5 ),  & ( 3, -12, 5, 8, -3 ),  & ( 3, -12, 5, 10, -5 ), \\
( 3, -12, 6, 7, -3 ),  & ( 3, -12, 6, 9, -5 ),  & ( 3, -12, 7, 6, -3 ),  & ( 3, -12, 7, 8, -5 ), \\
( 3, -9, -3, 10, 0 ),  & ( 3, -9, -3, 12, -2 ),  & ( 3, -9, -2, 9, 0 ),  & ( 3, -9, -2, 11, -2 ), \\
( 3, -9, -2, 13, -4 ),  & ( 3, -9, -2, 15, -6 ),  & ( 3, -9, -1, 8, 0 ),  & ( 3, -9, -1, 10, -2 ), \\
( 3, -9, -1, 12, -4 ),  & ( 3, -9, -1, 14, -6 ),  & ( 3, -9, 0, 7, 0 ),  & ( 3, -9, 0, 9, -2 ), \\
( 3, -9, 0, 11, -4 ),  & ( 3, -9, 0, 13, -6 ),  & ( 3, -9, 1, 6, 0 ),  & ( 3, -9, 1, 8, -2 ), \\
( 3, -9, 1, 10, -4 ),  & ( 3, -9, 1, 12, -6 ),  & ( 3, -9, 2, 5, 0 ),  & ( 3, -9, 2, 7, -2 ), \\
( 3, -9, 2, 9, -4 ),  & ( 3, -9, 2, 11, -6 ),  & ( 3, -9, 3, 4, 0 ),  & ( 3, -9, 3, 6, -2 ), \\
( 3, -9, 3, 8, -4 ),  & ( 3, -9, 3, 10, -6 ),  & ( 3, -9, 4, 3, 0 ),  & ( 3, -9, 4, 5, -2 ), \\
\end{array}}
$$

$$
{\tiny
\begin{array}{llll}
( 3, -9, 4, 7, -4 ),  & ( 3, -9, 4, 9, -6 ),  & ( 3, -9, 5, 2, 0 ),  & ( 3, -9, 5, 4, -2 ), \\
( 3, -9, 5, 6, -4 ),  & ( 3, -9, 5, 8, -6 ),  & ( 3, -9, 6, 3, -2 ),  & ( 3, -9, 6, 5, -4 ), \\
( 3, -9, 6, 7, -6 ),  & ( 3, -9, 7, 4, -4 ),  & ( 3, -9, 7, 6, -6 ),  & ( 3, -6, -6, 13, -3 ), \\
( 3, -6, -5, 8, 1 ),  & ( 3, -6, -5, 12, -3 ),  & ( 3, -6, -4, 7, 1 ),  & ( 3, -6, -4, 11, -3 ), \\
( 3, -6, -3, 4, 3 ),  & ( 3, -6, -3, 6, 1 ),  & ( 3, -6, -3, 10, -3 ),  & ( 3, -6, -2, 3, 3 ), \\
( 3, -6, -2, 5, 1 ),  & ( 3, -6, -2, 7, -1 ),  & ( 3, -6, -2, 9, -3 ),  & ( 3, -6, -1, 2, 3 ), \\
( 3, -6, -1, 4, 1 ),  & ( 3, -6, -1, 6, -1 ),  & ( 3, -6, -1, 8, -3 ),  & ( 3, -6, 0, 1, 3 ), \\
( 3, -6, 0, 3, 1 ),  & ( 3, -6, 0, 5, -1 ),  & ( 3, -6, 0, 7, -3 ),  & ( 3, -6, 1, 0, 3 ), \\
( 3, -6, 1, 2, 1 ),  & ( 3, -6, 1, 4, -1 ),  & ( 3, -6, 1, 6, -3 ),  & ( 3, -6, 2, 1, 1 ), \\
( 3, -6, 2, 3, -1 ),  & ( 3, -6, 2, 5, -3 ),  & ( 3, -6, 3, 0, 1 ),  & ( 3, -6, 3, 2, -1 ), \\
( 3, -6, 3, 4, -3 ),  & ( 3, -6, 4, 3, -3 ),  & ( 3, -6, 5, 2, -3 ),  & ( 3, -3, -6, 7, 0 ), \\
( 3, -3, -5, 2, 4 ),  & ( 3, -3, -5, 6, 0 ),  & ( 3, -3, -4, 1, 4 ),  & ( 3, -3, -4, 5, 0 ), \\
( 3, -3, -3, 0, 4 ),  & ( 3, -3, -3, 4, 0 ),  & ( 3, -3, -2, -1, 4 ),  & ( 3, -3, -2, 3, 0 ), \\
( 3, -3, -1, -2, 4 ),  & ( 3, -3, -1, 2, 0 ),  & ( 3, -3, 0, 1, 0 ),  & ( 3, -3, 1, 0, 0 ), \\
( 4, -12, 0, 10, -1 ),  & ( 4, -12, 1, 9, -1 ),  & ( 4, -12, 1, 11, -3 ),  & ( 4, -12, 2, 8, -1 ), \\
( 4, -12, 2, 10, -3 ),  & ( 4, -12, 3, 7, -1 ),  & ( 4, -12, 3, 9, -3 ),  & ( 4, -12, 4, 6, -1 ), \\
( 4, -12, 4, 8, -3 ),  & ( 4, -12, 5, 5, -1 ),  & ( 4, -12, 5, 7, -3 ),  & ( 4, -9, -3, 13, -4 ), \\
( 4, -9, -2, 12, -4 ),  & ( 4, -9, -1, 11, -4 ),  & ( 4, -9, 0, 6, 0 ),  & ( 4, -9, 0, 10, -4 ), \\
( 4, -9, 1, 5, 0 ),  & ( 4, -9, 1, 9, -4 ),  & ( 4, -9, 2, 4, 0 ),  & ( 4, -9, 2, 8, -4 ), \\
( 4, -9, 3, 7, -4 ),  & ( 4, -9, 4, 6, -4 ),  & ( 4, -9, 5, 5, -4 ),  & ( 4, -6, -3, 7, -1 ), \\
( 4, -6, -2, 2, 3 ),  & ( 4, -6, -2, 6, -1 ),  & ( 4, -6, -1, 1, 3 ),  & ( 4, -6, -1, 5, -1 ), \\
( 4, -6, 0, 4, -1 ),  & ( 4, -6, 1, 3, -1 ),  \\
\end{array}}
$$
 
\end{document}